\documentclass[12pt]{article}

\usepackage{CJK}

\usepackage{amsfonts,amssymb,amsthm,amsmath,latexsym}
\usepackage{subeqnarray,indent,eepicemu,url,cite,bm}
\usepackage{multirow}  
\usepackage{tensor}  
\usepackage{float}  

\date{April 12, 2018 \\[2mm]
      To appear in the {\em American Mathematical Monthly}\/}

\begin{document}

\title{How to generalize (and not to generalize) \\
       the Chu--Vandermonde identity}

\author{
   {\small Alan D.~Sokal}                  \\[2mm]
   {\small\it Department of Physics}       \\[-2mm]
   {\small\it New York University}         \\[-2mm]
   {\small\it 4 Washington Place}          \\[-2mm]
   {\small\it New York, NY 10003}          \\[-2mm]
   {\small\it USA}                         \\[-2mm]
   {\small\tt sokal@nyu.edu}               \\[-2mm]
   {\protect\makebox[5in]{\quad}}  
   \\[-2mm]
   {\small\it Department of Mathematics}   \\[-2mm]
   {\small\it University College London}   \\[-2mm]
   {\small\it Gower Street}                \\[-2mm]
   {\small\it London WC1E 6BT}             \\[-2mm]
   {\small\it UNITED KINGDOM}              \\[-2mm]
   {\small\tt sokal@math.ucl.ac.uk}        \\[3mm]
}

\maketitle
\thispagestyle{empty}   

\vspace*{-5mm}

\begin{abstract}
We consider two different interpretations of the
Chu--Vandermonde identity:
as an identity for polynomials, and as an identity for infinite matrices.
Each interpretation leads to a class of possible generalizations,
and in both cases we obtain a complete characterization of the solutions.
\end{abstract}


\bigskip
\bigskip
\noindent
{\bf Mathematics Subject Classification (MSC 2010) codes:}
05A19 (Primary);
05A10, 05A15, 05A40, 13B25, 15A24, 15B36, 33C47 (Secondary).

\clearpage

\newtheorem{theorem}{Theorem}
\newtheorem{proposition}[theorem]{Proposition}
\newtheorem{lemma}[theorem]{Lemma}
\newtheorem{corollary}[theorem]{Corollary}
\newtheorem{definition}[theorem]{Definition}
\newtheorem{conjecture}[theorem]{Conjecture}
\newtheorem{question}[theorem]{Question}
\newtheorem{problem}[theorem]{Problem}
\newtheorem{example}[theorem]{Example}

\renewcommand{\theenumi}{\alph{enumi}}
\renewcommand{\labelenumi}{(\theenumi)}
\def\eop{\hbox{\kern1pt\vrule height6pt width4pt
depth1pt\kern1pt}\medskip}
\def\prf{\par\noindent{\bf Proof.\enspace}\rm}
\def\rmk{\par\medskip\noindent{\bf Remark\enspace}\rm}

\newcommand{\textbfit}[1]{\textbf{\textit{#1}}}

\newcommand{\bigdash}{%
\smallskip\begin{center} \rule{5cm}{0.1mm} \end{center}\smallskip}

\newcommand{\safepar}{ {\protect\hfill\protect\break\hspace*{5mm}} }

\newcommand{\be}{\begin{equation}}
\newcommand{\ee}{\end{equation}}
\newcommand{\<}{\langle}
\renewcommand{\>}{\rangle}
\newcommand{\widebar}{\overline}
\def\reff#1{(\protect\ref{#1})}
\def\spose#1{\hbox to 0pt{#1\hss}}
\def\ltapprox{\mathrel{\spose{\lower 3pt\hbox{$\mathchar"218$}}
    \raise 2.0pt\hbox{$\mathchar"13C$}}}
\def\gtapprox{\mathrel{\spose{\lower 3pt\hbox{$\mathchar"218$}}
    \raise 2.0pt\hbox{$\mathchar"13E$}}}
\def\textprime{${}^\prime$}
\def\proof{\par\medskip\noindent{\sc Proof.\ }}
\def\firstproof{\par\medskip\noindent{\sc First Proof.\ }}
\def\secondproof{\par\medskip\noindent{\sc Second Proof.\ }}
\def\alternateproof{\par\medskip\noindent{\sc Alternate Proof.\ }}
\def\algebraicproof{\par\medskip\noindent{\sc Algebraic Proof.\ }}
\def\combinatorialproof{\par\medskip\noindent{\sc Combinatorial Proof.\ }}
\def\proofof#1{\bigskip\noindent{\sc Proof of #1.\ }}
\def\firstproofof#1{\bigskip\noindent{\sc First Proof of #1.\ }}
\def\secondproofof#1{\bigskip\noindent{\sc Second Proof of #1.\ }}
\def\thirdproofof#1{\bigskip\noindent{\sc Third Proof of #1.\ }}
\def\algebraicproofof#1{\bigskip\noindent{\sc Algebraic Proof of #1.\ }}
\def\combinatorialproofof#1{\bigskip\noindent{\sc Combinatorial Proof of #1.\ }}
\def\sketchofproof{\par\medskip\noindent{\sc Sketch of proof.\ }}
\renewcommand{\qed}{ $\square$ \bigskip}
\newcommand{\myendremark}{ $\blacksquare$ \bigskip}
\def\half{ {1 \over 2} }
\def\third{ {1 \over 3} }
\def\twothird{ {2 \over 3} }
\def\smfrac#1#2{{\textstyle{#1\over #2}}}
\def\smhalf{ {\smfrac{1}{2}} }
\def\smquarter{ {\smfrac{1}{4}} }
\newcommand{\real}{\mathop{\rm Re}\nolimits}
\renewcommand{\Re}{\mathop{\rm Re}\nolimits}
\newcommand{\imag}{\mathop{\rm Im}\nolimits}
\renewcommand{\Im}{\mathop{\rm Im}\nolimits}
\newcommand{\sgn}{\mathop{\rm sgn}\nolimits}
\newcommand{\tr}{\mathop{\rm tr}\nolimits}
\newcommand{\tg}{\mathop{\rm tg}\nolimits}
\newcommand{\supp}{\mathop{\rm supp}\nolimits}
\newcommand{\disc}{\mathop{\rm disc}\nolimits}
\newcommand{\diag}{\mathop{\rm diag}\nolimits}
\newcommand{\csch}{\mathop{\rm csch}\nolimits}
\newcommand{\tridiag}{\mathop{\rm tridiag}\nolimits}
\newcommand{\AZ}{\mathop{\rm AZ}\nolimits}
\newcommand{\perm}{\mathop{\rm perm}\nolimits}
\def\hboxscript#1{ {\hbox{\scriptsize\em #1}} }
\renewcommand{\emptyset}{\varnothing}
\newcommand{\eqdef}{\stackrel{\rm def}{=}}

\newcommand{\restrict}{\upharpoonright}

\newcommand{\compinv}{{\langle -1 \rangle}}   

\newcommand{\scra}{{\mathcal{A}}}
\newcommand{\scrb}{{\mathcal{B}}}
\newcommand{\scrc}{{\mathcal{C}}}
\newcommand{\scrd}{{\mathcal{D}}}
\newcommand{\scre}{{\mathcal{E}}}
\newcommand{\scrf}{{\mathcal{F}}}
\newcommand{\scrg}{{\mathcal{G}}}
\newcommand{\scrh}{{\mathcal{H}}}
\newcommand{\scri}{{\mathcal{I}}}
\newcommand{\scrj}{{\mathcal{J}}}
\newcommand{\scrk}{{\mathcal{K}}}
\newcommand{\scrl}{{\mathcal{L}}}
\newcommand{\scrm}{{\mathcal{M}}}
\newcommand{\scrn}{{\mathcal{N}}}
\newcommand{\scro}{{\mathcal{O}}}
\newcommand{\scrp}{{\mathcal{P}}}
\newcommand{\scrq}{{\mathcal{Q}}}
\newcommand{\scrr}{{\mathcal{R}}}
\newcommand{\scrs}{{\mathcal{S}}}
\newcommand{\scrt}{{\mathcal{T}}}
\newcommand{\scrv}{{\mathcal{V}}}
\newcommand{\scrw}{{\mathcal{W}}}
\newcommand{\scrz}{{\mathcal{Z}}}

\newcommand{\bfa}{{\mathbf{a}}}
\newcommand{\bfb}{{\mathbf{b}}}
\newcommand{\bfc}{{\mathbf{c}}}
\newcommand{\bfd}{{\mathbf{d}}}
\newcommand{\bfe}{{\mathbf{e}}}
\newcommand{\bfj}{{\mathbf{j}}}
\newcommand{\bfi}{{\mathbf{i}}}
\newcommand{\bfk}{{\mathbf{k}}}
\newcommand{\bfl}{{\mathbf{l}}}
\newcommand{\bfm}{{\mathbf{m}}}
\newcommand{\bfx}{{\mathbf{x}}}
\renewcommand{\k}{{\mathbf{k}}}
\newcommand{\n}{{\mathbf{n}}}
\newcommand{\vv}{{\mathbf{v}}}
\newcommand{\bv}{{\mathbf{v}}}
\newcommand{\w}{{\mathbf{w}}}
\newcommand{\x}{{\mathbf{x}}}
\newcommand{\cc}{{\mathbf{c}}}
\newcommand{\zero}{{\mathbf{0}}}
\newcommand{\one}{{\mathbf{1}}}
\newcommand{\bmm}{{\mathbf{m}}}

\newcommand{\ahat}{{\widehat{a}}}
\newcommand{\Zhat}{{\widehat{Z}}}

\newcommand{\C}{{\mathbb C}}
\newcommand{\D}{{\mathbb D}}
\newcommand{\Z}{{\mathbb Z}}
\newcommand{\N}{{\mathbb N}}
\newcommand{\Q}{{\mathbb Q}}
\newcommand{\PP}{{\mathbb P}}
\newcommand{\R}{{\mathbb R}}
\newcommand{\RR}{{\mathbb R}}
\newcommand{\E}{{\mathbb E}}

\newcommand{\ba}{{\bm{a}}}
\newcommand{\bahat}{{\widehat{\bm{a}}}}
\newcommand{\sfa}{{{\sf a}}}
\newcommand{\bb}{{\bm{b}}}
\newcommand{\bc}{{\bm{c}}}
\newcommand{\bchat}{{\widehat{\bm{c}}}}
\newcommand{\bd}{{\bm{d}}}
\newcommand{\bee}{{\bm{e}}}
\newcommand{\bff}{{\bm{f}}}
\newcommand{\bg}{{\bm{g}}}
\newcommand{\bh}{{\bm{h}}}
\newcommand{\bll}{{\bm{\ell}}}
\newcommand{\bp}{{\bm{p}}}
\newcommand{\br}{{\bm{r}}}
\newcommand{\bs}{{\bm{s}}}
\newcommand{\bu}{{\bm{u}}}
\newcommand{\bw}{{\bm{w}}}
\newcommand{\bx}{{\bm{x}}}
\newcommand{\by}{{\bm{y}}}
\newcommand{\bz}{{\bm{z}}}
\newcommand{\bA}{{\bm{A}}}
\newcommand{\bB}{{\bm{B}}}
\newcommand{\bC}{{\bm{C}}}
\newcommand{\bE}{{\bm{E}}}
\newcommand{\bF}{{\bm{F}}}
\newcommand{\bG}{{\bm{G}}}
\newcommand{\bH}{{\bm{H}}}
\newcommand{\bI}{{\bm{I}}}
\newcommand{\bJ}{{\bm{J}}}
\newcommand{\bM}{{\bm{M}}}
\newcommand{\bN}{{\bm{N}}}
\newcommand{\bP}{{\bm{P}}}
\newcommand{\bQ}{{\bm{Q}}}
\newcommand{\bS}{{\bm{S}}}
\newcommand{\bT}{{\bm{T}}}
\newcommand{\bW}{{\bm{W}}}
\newcommand{\bX}{{\bm{X}}}
\newcommand{\bIB}{{\bm{IB}}}
\newcommand{\bOB}{{\bm{OB}}}
\newcommand{\bOS}{{\bm{OS}}}
\newcommand{\bERR}{{\bm{ERR}}}
\newcommand{\bSP}{{\bm{SP}}}
\newcommand{\bMV}{{\bm{MV}}}
\newcommand{\bBM}{{\bm{BM}}}
\newcommand{\balpha}{{\bm{\alpha}}}
\newcommand{\bbeta}{{\bm{\beta}}}
\newcommand{\bgamma}{{\bm{\gamma}}}
\newcommand{\bdelta}{{\bm{\delta}}}
\newcommand{\bkappa}{{\bm{\kappa}}}
\newcommand{\bomega}{{\bm{\omega}}}
\newcommand{\bsigma}{{\bm{\sigma}}}
\newcommand{\btau}{{\bm{\tau}}}
\newcommand{\bpsi}{{\bm{\psi}}}
\newcommand{\bzeta}{{\bm{\zeta}}}
\newcommand{\bone}{{\bm{1}}}
\newcommand{\bzero}{{\bm{0}}}

%
%
\newcommand{\stirlingsubset}[2]{\genfrac{\{}{\}}{0pt}{}{#1}{#2}}
\newcommand{\stirlingcycleold}[2]{\genfrac{[}{]}{0pt}{}{#1}{#2}}
\newcommand{\stirlingcycle}[2]{\left[\! \stirlingcycleold{#1}{#2} \!\right]}
\newcommand{\assocstirlingsubset}[3]{{\genfrac{\{}{\}}{0pt}{}{#1}{#2}}_{\! \ge #3}}
\newcommand{\genstirlingsubset}[4]{{\genfrac{\{}{\}}{0pt}{}{#1}{#2}}_{\! #3,#4}}
\newcommand{\euler}[2]{\genfrac{\langle}{\rangle}{0pt}{}{#1}{#2}}
\newcommand{\eulergen}[3]{{\genfrac{\langle}{\rangle}{0pt}{}{#1}{#2}}_{\! #3}}
\newcommand{\eulersecond}[2]{\left\langle\!\! \euler{#1}{#2} \!\!\right\rangle}
\newcommand{\eulersecondgen}[3]{{\left\langle\!\! \euler{#1}{#2} \!\!\right\rangle}_{\! #3}}
\newcommand{\binomvert}[2]{\genfrac{\vert}{\vert}{0pt}{}{#1}{#2}}
\newcommand{\binomsquare}[2]{\genfrac{[}{]}{0pt}{}{#1}{#2}}


\newenvironment{sarray}{
             \textfont0=\scriptfont0
             \scriptfont0=\scriptscriptfont0
             \textfont1=\scriptfont1
             \scriptfont1=\scriptscriptfont1
             \textfont2=\scriptfont2
             \scriptfont2=\scriptscriptfont2
             \textfont3=\scriptfont3
             \scriptfont3=\scriptscriptfont3
           \renewcommand{\arraystretch}{0.7}
           \begin{array}{l}}{\end{array}}

\newenvironment{scarray}{
             \textfont0=\scriptfont0
             \scriptfont0=\scriptscriptfont0
             \textfont1=\scriptfont1
             \scriptfont1=\scriptscriptfont1
             \textfont2=\scriptfont2
             \scriptfont2=\scriptscriptfont2
             \textfont3=\scriptfont3
             \scriptfont3=\scriptscriptfont3
           \renewcommand{\arraystretch}{0.7}
           \begin{array}{c}}{\end{array}}

\clearpage

\section*{1. Introduction.}

One of the most celebrated formulae of elementary combinatorics
is the Chu--Vandermonde identity
\be
   \sum_{k=0}^n \binom{x}{k} \binom{y}{n-k}  \;=\;  \binom{x+y}{n}
 \label{eq.vandermonde}
\ee
(where $n$ is a nonnegative integer);%
\begin{CJK}{UTF8}{bsmi}
it was found by Chu Shih-Chieh (Zh\={u} Sh\`{\i}ji\'e, 朱世杰)%
\end{CJK}
%
sometime before 1303 \cite{Zhu_2006}
and was rediscovered circa 1772 by Alexandre-Th\'eophile Vandermonde
\cite{Vandermonde_1772}.
If we take $x$ and $y$ to be nonnegative integers,
this identity has a well-known combinatorial interpretation and proof:
we can choose an $n$-person committee from a group of $x$ women and $y$ men
in $\binom{x+y}{n}$ ways;
on the other hand, if this committee is to have $k$ women and $n-k$ men,
then the number of ways is $\binom{x}{k} \binom{y}{n-k}$;
summing over all possible values of $k$ gives the result.
But now, having proven the identity for nonnegative integer $x$~and~$y$,
we can shift gears and regard $x$~and~$y$ as algebraic indeterminates:
then both sides of \reff{eq.vandermonde} are polynomials in $x$~and~$y$
(of total degree~$n$),
which agree whenever $x$~and~$y$ are nonnegative integers;
since two polynomials that agree at infinitely many points must be equal,
it follows that \reff{eq.vandermonde} holds as a polynomial identity
in $x$~and~$y$.

We thus see already from this brief account that the Chu--Vandermonde identity
can be given at least two distinct interpretations.
On the one hand, we can regard $x$~and~$y$ as indeterminates;
then \reff{eq.vandermonde} is the identity
\be
    \sum_{k=0}^n f_k(x) \, f_{n-k}(y)  \;=\;  f_n(x+y)
 \label{eq.fksum}
\ee
for the polynomials
\be
   f_n(x)  \;=\;  \binom{x}{n}  \;\eqdef\; {x^{\underline{n}} \over n!}
           \;\eqdef\; {x(x-1) \cdots (x-n+1) \over n!}
           \;\in\; \Q[x]
   \;.
 \label{def.f.binomial}
\ee
(Of course, $x$ and $y$ can then be specialized, if we wish,
to specific values in any commutative ring containing the rationals,
such as the real or complex numbers.)
On the other hand, we can restrict $x$ and $y$ to be nonnegative integers;
then \reff{eq.vandermonde} is the identity
\be
   \sum_{j=0}^n L_{ij} \, L_{\ell,n-j}  \;=\;  L_{i+\ell,n}
 \label{eq.identity.L}
\ee
for the infinite lower-triangular Pascal matrix
\cite{Call_93,Aceto_01,Edelman_04}
$L = (L_{ij})_{i,j \ge 0}$ defined by
\be
   L_{ij}  \;=\;  \binom{i}{j}
   \;.
 \label{def.L.pascal}
\ee

At this point it is natural to ask for generalizations of
\reff{eq.fksum}/\reff{def.f.binomial}
and \reff{eq.identity.L}/\reff{def.L.pascal}, respectively.
Are there other sequences $\bff = (f_n)_{n \ge 0}$
of polynomials (or formal power series) satisfying \reff{eq.fksum},
and if so, can we classify them all?
Likewise, are there other lower-triangular matrices $L = (L_{ij})_{i,j \ge 0}$
satisfying \reff{eq.identity.L}, and if so, can we classify them?
The answer to the first question is well known,
and we will review it here (and extend it slightly).
Then we will give a negative answer to the second question,
but with an interesting twist;
some of these latter results seem to be new.

\bigskip

\section*{2. Convolution families.}

Let $R$ be a commutative ring,
and let $\bff = (f_n)_{n \ge 0}$
be a sequence of formal power series $f_n(x) \in R[[x]]$.
We will call $\bff$ a {\em convolution family}\/
\cite{Knuth_92} (see also \cite{Labelle_80,Zeng_96})
if it satisfies \reff{eq.fksum};
and we will call $\bff$ a {\em weak convolution family}\/
if it satisfies \reff{eq.fksum} restricted to $x=y$:
\be
    \sum_{k=0}^n f_k(x) \, f_{n-k}(x)  \;=\;  f_n(2x)
   \;.
 \label{eq.fksum.weak}
\ee
We use the notation $[t^n] \, F(t)$
to denote the coefficient of $t^n$ in the formal power series (or polynomial)
$F(t) \in R[[t]]$.
We then have the following characterization of (weak) convolution families:

\begin{proposition}[Characterization of convolution families]
   \label{prop.convfamily}
Let $R$ be a commutative ring containing the rationals,
and let $\Psi(t) = \sum_{n=0}^\infty \psi_n t^n \in R[[t]]$
be any formal power series.
Then
\be
   f_n(x)  \;=\;  [t^n] \, e^{x \Psi(t)}
 \label{eq.fn.G.H}
\ee
defines a convolution family in $R[[x]]$.
Moreover, $\widetilde{f}_n(x) \eqdef f_n(x) / e^{\psi_0 x}$
are polynomials in $x$ with $\deg \widetilde{f}_n \le n$
and satisfying $\widetilde{f}_n(0) = 0$ for $n \ge 1$.
[In~particular, the $f_n$ are themselves polynomials for all $n$
 if and only if $\psi_0$ is nilpotent,
 and are polynomials of degree at most $n$ for all $n$
 if and only if $\psi_0 = 0$.]

Conversely, let $R$ be a commutative ring containing the rationals
in which the only idempotents are 0 and 1
(for instance, an integral domain or a local ring),
and let $\bff = (f_n)_{n \ge 0}$ be a weak convolution family in $R[[x]]$.
Then either $\bff$ is identically zero,
or else there exists a unique formal power series $\Psi \in R[[t]]$
such that \reff{eq.fn.G.H} holds.
In~particular, $\bff$ is actually a convolution family,
and the foregoing statements about $\widetilde{f}_n$ hold.
\end{proposition}

\proof
By taking the coefficient of $t^n$ on both sides of the identity
$e^{x \Psi(t)} e^{y \Psi(t)} = e^{(x+y) \Psi(t)}$,
we see, using \reff{eq.fn.G.H}, that
\be
    \sum_{k=0}^n f_k(x) \, f_{n-k}(y)  \;=\;  f_n(x+y)
    \;,
\ee
or in other words that $\bff$ is a convolution family.
Moreover, it is easy to see that
$f_n(x)/e^{\psi_0 x} = [t^n] \, e^{x [\Psi(t) - \psi_0]}$
is a polynomial in $x$ of degree at most $n$,
which has zero constant term for $n \ge 1$.
For the statements in brackets, ``if'' is easy,
and ``only if'' follows by looking at $n=0$.

%

For the converse, let us use the notation $f_{nk} = [x^k] \, f_n(x)$.
By using \reff{eq.fksum.weak} specialized to $n=0$ and $x=0$,
we see that $f_{00}^2 = f_{00}$;
so the hypothesis on $R$ implies that either $f_{00} = 0$ or $f_{00} = 1$.
If $f_{00} = 0$, then a simple inductive argument using
\reff{eq.fksum.weak} specialized to $n=0$
shows that $f_{0k} = 0$ for all $k$, i.e., $f_0 = 0$;
and then a second inductive argument using \reff{eq.fksum.weak}
shows that $f_n = 0$ for all $n$.
So we can assume that $f_{00} = 1$.
Define the bivariate ordinary generating function
$F(x,t) = \sum_{n=0}^\infty f_n(x) \, t^n \in R[[x,t]]$.
Since $F$ has constant term $f_{00} = 1$,
we can define $L(x,t) = \log F(x,t) \in R[[x,t]]$,
which has constant term 0.
The identity \reff{eq.fksum.weak}, translated to generating functions,
says that $F(x,t)^2 = F(2x,t)$, or equivalently that $2 L(x,t) = L(2x,t)$.
Writing $L(x,t) = \sum_{k=0}^\infty \ell_k(t) \, x^k$,
we see that $\sum_{k=0}^\infty (2 - 2^k) \, \ell_k(t) \, x^k = 0$
as a power series in $x$,
hence that $\ell_k(t) = 0$ for all $k \neq 1$.
It follows that $L(x,t) = x \Psi(t)$ for some $\Psi \in R[[t]]$.
Therefore $f_n(x) = [t^n] \, e^{x \Psi(t)}$,
as claimed in \reff{eq.fn.G.H}.
Indeed, \reff{eq.fn.G.H} says precisely that $L(x,t) = x \Psi(t)$,
so $\Psi$ is uniquely determined by $\bff$.
\qed

\bigskip

\noindent
{\bf Examples of convolution families.}
1.  The binomial theorem corresponds to $f_n(x) = x^n/n!$
and $\Psi(t) = t$.

2.  The Chu--Vandermonde identity corresponds to
$f_n(x) = \binom{x}{n} = x^{\underline{n}}/n!$
and $\Psi(t) = \log(1 + t)$.

3.  The ``dual Chu--Vandermonde identity'' corresponds to
$f_n(x) = \binom{x+n-1}{n} = x^{\overline{n}}/n!$
and $\Psi(t) = -\log(1 - t)$.

4.  Define the univariate Bell polynomials
$B_n(x) = \sum_{k=0}^n \stirlingsubset{n}{k} x^k$,
where the Stirling number $\stirlingsubset{n}{k}$
is the number of partitions of an $n$-element set into $k$ nonempty blocks
\cite{Graham_94}.
Then $f_n(x) = B_n(x)/n!$ corresponds to $\Psi(t) = e^t - 1$.

\bigskip

\noindent
See \cite{Knuth_92} for many further examples.

\bigskip
\bigskip

\noindent
{\bf Remarks.}
1.  Proposition~\ref{prop.convfamily} is a slight extension of
\cite[Theorem~4.1]{Garsia_73}, \cite{Knuth_92}
and \cite[Exercise~5.37a]{Stanley_99},
in two directions:
allowing the $f_n$ to be formal power series rather than just polynomials,
and allowing the coefficients to lie in a commutative ring containing
the rationals rather than just a field of characteristic~0.
Knuth \cite{Knuth_92} proved the converse half by an inductive argument
based on carefully examining the coefficients $f_{nk}$.
This is nice, but it seems to me that the generating-function argument
using the logarithm, taken from \cite{Garsia_73,Stanley_99},
is simpler and more enlightening.

2.  For the converse, it would alternatively suffice to assume that
\reff{eq.fksum} holds for $y = rx$
for {\em any}\/ fixed rational number $r \neq 0,-1$ (not just $r=1$):
we would then have $L(x,t) + L(rx,t) = L((1+r)x,t)$,
which also implies $L(x,t) = x \Psi(t)$
since $1 + r^m \neq (1+r)^m$ for all $m \in \N \setminus \{1\}$
and all $r \in \Q \setminus \{0,-1\}$.

3.  By iterating \reff{eq.fksum}, we see that every convolution family $\bff$
satisfies
\be
    \sum\limits_{\begin{scarray}
                    k_1,\ldots,k_m \ge 0 \\
                    k_1 + \ldots + k_m = n
                 \end{scarray}}
    \!\!\!\!\!
    f_{k_1}(x_1) \,\cdots\, f_{k_m}(x_m)  \;=\;  f_n(x_1 + \ldots + x_m)
 \label{eq.fksum.iterated}
\ee
for each integer $m \ge 1$.
This ``multinomial'' version of the identity is sometimes useful.

4.  The foregoing theory can equivalently be expressed
in terms of the polynomials $F_n(x) = n! \, f_n(x)$, which satisfy
\be
    \sum_{k=0}^n \binom{n}{k} F_k(x) \, F_{n-k}(y)  \;=\;  F_n(x+y)
    \;.
 \label{eq.Fksum.binomial}
\ee
Sequences of polynomials satisfying \reff{eq.Fksum.binomial}
are called {\em sequences of binomial type}\/;
they were introduced by Rota and collaborators
\cite{Mullin_70,Rota_73,Garsia_73,Fillmore_73,Roman_78,Roman_84}
and studied by means of the umbral calculus
\cite{Roman_78,Roman_84,Rota_94,Gessel_03}.\footnote{
   As part of the definition of ``sequence of binomial type'',
   Rota {\em et al.}\/ \cite{Mullin_70,Rota_73}
   and many subsequent authors \cite{Garsia_73,Fillmore_73,Roman_78,Roman_84}
   imposed the additional condition that $\deg F_n = n$ exactly.
   But this condition is irrelevant for our purposes,
   so we refrain from imposing it.
}
A purely combinatorial approach to sequences of binomial type,
employing the theory of species,
has been developed by Labelle \cite{Labelle_81}
(see also \cite[Section~3.1]{Bergeron_98}).
See also \cite{Zeng_96,Scott-Sokal_expidentities}
for some multivariate generalizations.

The identity \reff{eq.Fksum.binomial}
is closely related to the {\em exponential formula}\/
\cite{Wilf_94,Bergeron_98,Stanley_99}
in enumerative combinatorics,
which relates the weights $c_n$ of ``connected'' objects
to the weights $F_n(x)$ of ``all'' objects,
when the weight of an object is defined to be the product of the weights
of its ``connected components'' times a factor $x$ for each
connected component.
See \cite{Wilf_94,Bergeron_98,Stanley_99,Scott-Sokal_expidentities}
for further discussion and many applications.
\myendremark

\bigskip

\noindent
{\bf Open questions.}
What happens if $R$ does not contain the rationals,
or if $R$ has idempotents other than 0 and 1?

\bigskip

\section*{3. Convolution families, generalized.}

Let us now generalize the identity \reff{eq.fksum}
by allowing three different sequences $f,g,h$ in place of $f,f,f$:
\be
   \sum_{k=0}^n f_k(x) \, g_{n-k}(y)  \;=\;  h_n(x+y)
   \;.
 \label{eq.conv.fgh}
\ee
Do any new solutions arise?
It turns out that the answer is yes, as follows:

\begin{proposition}
   \label{prop.convfamily.bis}
Let $R$ be a commutative ring containing the rationals,
and let $A,B,\Psi \in R[[t]]$.  Then
\be
   f_n(x)  \:=\:  [t^n] \, A(t) e^{x \Psi(t)}  \:,\quad
   g_n(x)  \:=\:  [t^n] \, B(t) e^{x \Psi(t)}  \:,\quad
   h_n(x)  \:=\:  [t^n] \, A(t) B(t) e^{x \Psi(t)}
 \label{eq.prop.convfamily.bis}
\ee
are sequences in $R[[x]]$ satisfying \reff{eq.conv.fgh}.
Moreover, $\widetilde{f}_n(x) \eqdef f_n(x) / e^{\psi_0 x}$
[where $\psi_0 \eqdef \Psi(0)$]
are polynomials in $x$ with $\deg \widetilde{f}_n \le n$,
and likewise for $g_n$ and $h_n$.
[In~particular, the $f_n,g_n,h_n$ are themselves polynomials for all $n$
 if and only if $\psi_0$ is nilpotent,
 and are polynomials of degree at most $n$ for all $n$
 if and only if $\psi_0 = 0$.]

Conversely, let $R$ be a commutative ring containing the rationals,
and let
$\bff = (f_n)_{n \ge 0}$, $\bg = (g_n)_{n \ge 0}$, $\bh = (h_n)_{n \ge 0}$
be sequences in $R[[x]]$ satisfying \reff{eq.conv.fgh}.
Assume further that $f_0(0)$ and $g_0(0)$ are invertible in $R$
[or equivalently that $h_0(0) = f_0(0) \, g_0(0)$ is invertible in $R$].
Then there exist unique formal power series $A,B,\Psi \in R[[t]]$
such that \reff{eq.prop.convfamily.bis} holds;
and $A$ and $B$ are invertible in $R[[t]]$.
\end{proposition}

\proof
By taking the coefficient of $t^n$ on both sides of the identity
$A(t) e^{x \Psi(t)} \, B(t) e^{y \Psi(t)} = {A(t) B(t) e^{(x+y) \Psi(t)}}$,
we see that \reff{eq.conv.fgh} holds.
The statements about $f_n,g_n,h_n$ are proven
as in Proposition~\ref{prop.convfamily}.

For the converse, let $\alpha = f_0(0)$ and $\beta = g_0(0)$;
by hypothesis $\alpha$ and $\beta$ are invertible, and $h_0(0) = \alpha\beta$.
Define the ordinary generating functions
$F(x,t) = \sum_{n=0}^\infty f_n(x) \, t^n$,
$G(x,t) = \sum_{n=0}^\infty g_n(x) \, t^n$
and $H(x,t) = \sum_{n=0}^\infty h_n(x) \, t^n$;
they have constant terms $\alpha$, $\beta$ and $\alpha\beta$, respectively.
Then define $L(x,t) = \log [\alpha^{-1} F(x,t)]$,
$M(x,t) = \log [\beta^{-1} G(x,t)]$ and
$N(x,t) = \log [(\alpha\beta)^{-1} H(x,t)]$;
they have constant term 0.
The identity \reff{eq.conv.fgh}, translated to generating functions,
says that $F(x,t) \, G(y,t) = H(x+y,t)$,
or equivalently that $L(x,t) + M(y,t) = N(x+y,t)$.
Then we must have $[x^k] N(x,t) = 0$ for $k \ge 2$
in order to avoid terms $x^i y^j$ with $i,j \ge 1$ in $N(x+y,t)$.
It then follows that we must have
\be
   L(x,t) \:=\: \Gamma(t) \,+\, x \Psi(t)  \:,\quad
   M(x,t) \:=\: \Delta(t) \,+\, x \Psi(t)  \:,\quad
   N(x,t) \:=\: \Gamma(t) \,+\, \Delta(t) \,+\, x \Psi(t)
\ee
for some formal power series $\Gamma,\Delta,\Psi \in R[[t]]$,
where $\Gamma$ and $\Delta$ have zero constant term.
Then \reff{eq.prop.convfamily.bis} holds with
$A(t) = \alpha e^{\Gamma(t)}$ and $B(t) = \beta e^{\Delta(t)}$.
Furthermore, $L,M,N$ are uniquely determined by $\bff,\bg,\bh$;
hence $\Gamma,\Delta,\Psi$ are also uniquely determined;
hence $A,B,\Psi$ are uniquely determined as well.
\qed

\bigskip

{\bf Remarks.}
1.  Under the assumption that $f_0(0)$ and $g_0(0)$ are invertible in $R$,
clearly $\bff = \bg = \bh$ if and only if $A(t) = B(t) = 1$;
and in this case we recover the situation of Proposition~\ref{prop.convfamily}.

2.  If the only idempotents in $R$ are 0 and 1,
then it is not difficult to show that the same holds true in $R[[t]]$. 
In this case $\bff = \bg = \bh$ if and only if
$A(t) = B(t) = 1$ or $A(t) = B(t) = 0$;
and we also recover the situation of Proposition~\ref{prop.convfamily}.
%
%

3.  Note that here, unlike in Proposition~\ref{prop.convfamily},
it does not suffice to assume \reff{eq.conv.fgh} only for $x=y$.
Indeed, the equality
\be
   \sum_{k=0}^n f_k(x) \, g_{n-k}(x)  \;=\;  h_n(2x)
\ee
can be taken as the {\em definition}\/ of $\bh$
for {\em completely arbitrary}\/ sequences $\bff$ and $\bg$.

4.  The foregoing theory can once again equivalently be expressed
in terms of the polynomials $F_n(x) = n! \, f_n(x)$
and likewise for $g_n$ and $h_n$, which satisfy
\be
    \sum_{k=0}^n \binom{n}{k} F_k(x) \, G_{n-k}(y)  \;=\;  H_n(x+y)
    \;.
 \label{eq.FGHsum.binomial}
\ee
Sequences of polynomials $F_n(x) = n! \, f_n(x)$
defined by \reff{eq.prop.convfamily.bis}
are called {\em Sheffer sequences}\/ \cite[pp.~2, 18--19]{Roman_84}.\footnote{
   As part of the definition of ``Sheffer sequence'',
   some authors (e.g.\ \cite[p.~2]{Roman_84})
   impose the additional conditions $A(0) \neq 0$,
   $B(0) = 0$, and $B'(0) \neq 0$.
   But these conditions are irrelevant for our purposes,
   so we refrain from imposing them.
}
Many classical sequences of polynomials ---
including the Hermite, Laguerre and Bernoulli polynomials ---
are Sheffer sequences \cite[pp.~2, 28--31, 53--130]{Roman_84}.
\myendremark

\bigskip

\noindent
{\bf Open question.}
What happens if $f_0(0)$ and/or $g_0(0)$ are not invertible?

\bigskip

\section*{4. Pascal-like matrices.}

We now turn to the matrix interpretation of the Chu--Vandermonde identity.
But before proceeding further,
let us generalize the identity \reff{eq.identity.L} in two ways:
First, we allow three different matrices $A,B,C$ in place of $L,L,L$;
and second, we do not require $A,B,C$ to be lower-triangular.

So let $A = (a_{ij})_{i,j \ge 0}$ be a matrix with entries
in a commutative ring $R$.  We define the {\em row-generating series}\/
$A_i(u) = \sum_{j=0}^\infty a_{ij} u^j \in R[[u]$.
Clearly, there is a one-to-one correspondence between matrices
$A = (a_{ij})_{i,j \ge 0}$ and collections $(A_i(u))_{i \ge 0}$
of row-generating series.

We then have the following result:

\begin{proposition}[Characterization of Pascal-like matrices]
   \label{prop.pascal-like}
Let $R$ be a commutative ring, and let $f,g,h \in R[[u]]$.
Then
\be
   A_i(u) \;=\; f(u) \, h(u)^i   \:,\quad
   B_i(u) \;=\; g(u) \, h(u)^i   \:,\quad
   C_i(u) \;=\; f(u) \, g(u) \, h(u)^i
 \label{eq.ABC.fgh}
\ee
are the row-generating series of matrices
$A = (a_{ij})_{i,j \ge 0}$,
$B = (b_{ij})_{i,j \ge 0}$,
$C = (c_{ij})_{i,j \ge 0}$
satisfying
\be
   \sum_{j=0}^n a_{ij} \, b_{\ell,n-j}  \;=\;  c_{i+\ell,n}
   \qquad\hbox{for all $i,\ell,n \ge 0$}
   \;.
 \label{eq.identity.ABC}
\ee

Conversely, suppose that
$A = (a_{ij})_{i,j \ge 0}$,
$B = (b_{ij})_{i,j \ge 0}$,
$C = (c_{ij})_{i,j \ge 0}$
are matrices with entries in a commutative ring $R$
that satisfy \reff{eq.identity.ABC}.
Suppose further that $a_{00}$ and $b_{00}$ are invertible in $R$
[or equivalently that $c_{00} = a_{00} b_{00}$ is invertible in $R$].
Then there exist unique series $f,g,h \in R[[u]]$
such that the row-generating series of $A,B,C$ are given by \reff{eq.ABC.fgh};
and $f$ and $g$ are invertible.
\end{proposition}

\proof
Multiplying \reff{eq.identity.ABC} by $u^n$ and summing over $n$,
we see that \reff{eq.identity.ABC} is equivalent to the equality
\be
   A_i(u) \, B_\ell(u)  \;=\;  C_{i+\ell}(u)
   \qquad\hbox{for all $i,\ell \ge 0$}
 \label{eq.AiBl}
\ee
for the row-generating series.
It is immediate that the construction \reff{eq.ABC.fgh}
satisfies \reff{eq.AiBl}.

For the converse, let us write out \reff{eq.AiBl} in detail:
\begin{subeqnarray}
   C_0  & = &  A_0 B_0  \\
   C_1  & = &  A_0 B_1  \;=\;  A_1 B_0  \\
   C_2  & = &  A_0 B_2  \;=\;  A_1 B_1  \;=\;  A_2 B_0  \\
        & \vdots &   \nonumber
\end{subeqnarray}
Since $a_{00}$ and $b_{00}$ are invertible in $R$,
it follows that $A_0$ and $B_0$ are invertible in $R[[u]]$.
From $C_n = A_0 B_n = A_n B_0$, we deduce that $A_n/A_0 = B_n/B_0$;
let us call this common value $h_n$.
Then \reff{eq.AiBl} says (after division by $A_0 B_0$) that
$h_i h_\ell = h_{i+\ell}$ for all $i,\ell$.
It follows by induction that $h_n = h_1^n$.
So \reff{eq.ABC.fgh} holds with $f = A_0$, $g = B_0$, $h = h_1$.
Moreover, it is clear from \reff{eq.ABC.fgh} that $f = A_0$ and $g = B_0$;
and since $A_0$ and $B_0$ are invertible, we must also have $h = A_1/A_0$.
So $f,g,h$ are uniquely determined.
\qed

\bigskip

\noindent
{\bf Open question.}
What happens if $a_{00}$ and/or $b_{00}$ are not invertible?

\bigskip
\bigskip

{\bf Remarks.}
1.  Under the assumption that $a_{00}$ and $b_{00}$ are invertible in $R$,
clearly $A=B=C$ if and only if $f(u) = g(u) = 1$.
It turns out (as I discovered after completing
the proof of Proposition~\ref{prop.pascal-like})
that the $A=B=C$ special case of Proposition~\ref{prop.pascal-like}
had been proven nearly 40 years ago
by Olive \cite[Theorems~3.2 and 3.3]{Olive_79}.

2.  If $f(t)$ and $h(t)$ are formal power series,
let $\scrr(f,h)$ be the infinite matrix
$(\scrr(f,h)_{nk})_{n,k \ge 0}$ defined by
\be
   \scrr(f,h)_{nk}
   \;=\;
   [t^n] \, f(t) h(t)^k
   \;.
 \label{def.riordan}
\ee
When $h$ has constant term 0, the matrix $\scrr(f,h)$ is lower-triangular
and is called a {\em Riordan array}\/ \cite{Shapiro_91,Sprugnoli_94,Barry_16};
such matrices arise frequently in enumerative combinatorics,
and the theory of Riordan arrays provides a useful unifying framework.
But --- as a handful of authors have noted
\cite{Bala_15} \cite{Hoggatt_75} \cite[p.~288]{Pemantle_13} ---
there are also interesting examples
in which $h$ has a nonzero constant term.
Let us call this more general concept,
in which $h$ is an arbitrary formal power series,
a {\em wide-sense Riordan array}\/.\footnote{
   References \cite{Pemantle_13,Bala_15}
   call this a ``generalized Riordan array'',
   but we prefer to avoid this term because it has already been used,
   in a highly-cited paper \cite{Wang_08},
   for a completely unrelated generalization of Riordan arrays.
}
We can then see that the matrices $A,B,C$ defined in \reff{eq.ABC.fgh}
are simply the transpose of a wide-sense Riordan array.

3. What is the relation between Propositions~\ref{prop.convfamily.bis}
and \ref{prop.pascal-like}?
Comparing \reff{eq.conv.fgh} with \reff{eq.identity.ABC},
we see that $x,y,k$ correspond to $i,\ell,j$, respectively;
and then, comparing \reff{eq.prop.convfamily.bis} with \reff{eq.ABC.fgh},
we see that $A(t), B(t), e^{\Psi(t)}$
correspond to $f(u), g(u), h(u)$, respectively.
Therefore, if $\Psi(0) = 0$, we can define $h = e^\Psi$ satisfying $h(0) = 1$;
and since in this case the $f_n,g_n,h_n$ are {\em polynomials}\/,
it makes sense to evaluate them at integer arguments
to obtain $f_n(i) = a_{in}$ and analogously for $g$ and $h$.
And conversely, if $h(0) = 1$,
we can define $\Psi = \log h$ satisfying $\Psi(0) = 0$,
and a triplet of matrices $A,B,C$ satisfying \reff{eq.ABC.fgh}
can be extended to a triplet $\bff,\bg,\bh$
of sequences of polynomials satisfying \reff{eq.prop.convfamily.bis};
once again we can evaluate them at integer arguments
to obtain $f_n(i) = a_{in}$ and analogously for $g$ and $h$.
In other words, when $h(0) = 1$, each column of the matrices $A,B,C$
is a {\em polynomial}\/ function of the index $i$.

If, by contrast, $\Psi(0) \eqdef \psi_0 \neq 0$ or $h(0) \neq 1$,
then Propositions~\ref{prop.convfamily.bis} and \ref{prop.pascal-like}
are incommensurable for general commutative rings $R$.
However, when $R = \R$ or $\C$, we can still define $h = e^\Psi$
satisfying $h(0) = e^{\psi_0}$;
then $f_n(x),g_n(x),h_n(x)$ are polynomials multiplied by $e^{\psi_0 x}$
and can again be evaluated at integer arguments.
And conversely, if $R= \R$ (resp.\ $\C$) and $h(0) > 0$ (resp.\ $h(0) \neq 0$),
then we can define $\Psi = \log h$,
and the columns of $A,B,C$ are interpolated by functions
$f_n(x), g_n(x), h_n(x)$ that are polynomials multiplied by $e^{\psi_0 x}$.
\myendremark

\bigskip

We can now specialize to the case in which $A$ and $B$ are lower-triangular.
In fact, we can be a bit more general.
Let us say that a matrix $M = (m_{ij})_{i,j \ge 0}$
is {\em lower-triangular in row $i$}\/ if $m_{ij} = 0$ for all $j > i$.
We then have:

\begin{corollary}
   \label{cor.pascal-like}
Let
$A = (a_{ij})_{i,j \ge 0}$,
$B = (b_{ij})_{i,j \ge 0}$,
$C = (c_{ij})_{i,j \ge 0}$
be matrices with entries in a commutative ring $R$
that satisfy \reff{eq.identity.ABC};
and suppose further that $a_{00}$ and $b_{00}$ are invertible in $R$.
If $A$ and $B$ are lower-triangular in row 0,
then there exist $\alpha,\beta \in R$ with $\alpha$ and $\beta$ invertible,
and $h \in R[[u]]$, such that
\be
   A_i(u) \;=\; \alpha \, h(u)^i   \:,\quad
   B_i(u) \;=\; \beta \, h(u)^i   \:,\quad
   C_i(u) \;=\; \alpha\beta \, h(u)^i
   \;.
\ee
If, in addition, at least one of $A,B,C$ is lower-triangular in row 1,
then there exist $\alpha,\beta,\kappa,\lambda \in R$,
with $\alpha$ and $\beta$ invertible, such that
\be
   a_{ij}  \;=\;  \alpha \kappa^{i-j} \lambda^j \binom{i}{j}  \:,\quad
   b_{ij}  \;=\;  \beta \kappa^{i-j} \lambda^j \binom{i}{j}  \:,\quad
   c_{ij}  \;=\;  \alpha\beta \kappa^{i-j} \lambda^j \binom{i}{j}
   \;.
\ee
\end{corollary}

\proof
If $A$ and $B$ are lower-triangular in row 0,
then $f = a_{00} = \alpha$ and $g = b_{00} = \beta$.
If, in addition, at least one of $A,B,C$ is lower-triangular in row 1,
then $h = \kappa + \lambda u$.
\qed

When $A=B=C$, this yields:

\begin{corollary}[No-go theorem for lower-triangular Pascal-like matrices]
   \label{cor.pascal-like.2}
Let $L = (L_{ij})_{i,j \ge 0}$ be a lower-triangular matrix
with entries in a commutative ring $R$
that satisfies \reff{eq.identity.L}.
If $L_{00}$ is invertible in $R$, then in fact $L_{00} = 1$,
and there exist $\kappa,\lambda \in R$ such that
\be
   L_{ij}  \;=\;  \kappa^{i-j} \lambda^j \binom{i}{j}
   \;.
\ee
\end{corollary}

Corollary~\ref{cor.pascal-like.2} thus gives a negative answer
to the question posed in the introduction:
there are no lower-triangular solutions to \reff{eq.identity.L}
other than trivial rescalings of the Pascal matrix.
But --- and this is the interesting twist ---
Proposition~\ref{prop.pascal-like} shows that there {\em are}\/
interesting examples if we give up the insistence that $L$ be lower-triangular:
we can take $f(u) = g(u) = 1$
and choose an {\em arbitrary}\/ formal power series $h(u)$,
not just $h(u) = \kappa + \lambda u$.
These examples turn out to have interesting applications
to the theory of Hankel-total positivity \cite{Sokal_totalpos};
but that story will have to be told elsewhere \cite{Petreolle-Sokal_nontri}.

\section*{Acknowledgments}

I wish to thank Mathias P\'etr\'eolle and Bao-Xuan Zhu
for helpful conversations,
and the referees for helpful suggestions.

This research was supported in part by
U.K.~Engineering and Physical Sciences Research Council grant EP/N025636/1.


\begin{thebibliography}{99}

\bibitem{Aceto_01}  L. Aceto and D. Trigiante,
   The matrices of Pascal and other greats,
   Amer. Math. Monthly {\bf 108}, 232--245 (2001).

\bibitem{Bala_15}  P. Bala, Notes on generalized Riordan arrays,
   \url{oeis.org/A260492/a260492.pdf}
   (August 2015).

\bibitem{Barry_16}  P. Barry, {\em Riordan Arrays: A Primer}\/
   (Logic Press, County Kildare, Ireland, 2016).

\bibitem{Bergeron_98}  F. Bergeron, G. Labelle and P. Leroux,
      {\em Combinatorial Species and Tree-Like Structures}\/
      (Cambridge University Press, Cambridge--New York, 1998).

\bibitem{Call_93}  G.~S. Call and D.~J. Velleman, Pascal's matrices,
   Amer. Math. Monthly {\bf 100}, 372--376 (1993).

\bibitem{Edelman_04}  A. Edelman and G. Strang, Pascal matrices,
   Amer. Math. Monthly {\bf 111}, 189--197 (2004).

\bibitem{Fillmore_73}  J.~P. Fillmore and S.~G. Williamson,
   A linear algebra setting for the Rota--Mullin theory of polynomials
   of binomial type,
   Linear and Multilinear Algebra {\bf 1}, 67--80 (1973).

\bibitem{Garsia_73}  A.~M. Garsia, An expos\'e of the Mullin--Rota theory
   of polynomials of binomial type,
   Linear and Multilinear Algebra {\bf 1}, 47--65 (1973).

\bibitem{Gessel_03}  I.~M. Gessel, Applications of the classical umbral
   calculus, Algebra Univers. {\bf 49}, 397--434 (2003).

\bibitem{Graham_94}  R.~L. Graham, D.~E. Knuth and O. Patashnik,
   {\em Concrete Mathematics: A Foundation for Computer Science}\/,
   2nd ed.~(Addison-Wesley, Reading, Mass., 1994).

\bibitem{Hoggatt_75}  V.~E. Hoggatt, Jr. and G.~E. Bergum,
   Generalized convolution arrays,
   Fibonacci Quart. {\bf 13}, 193--198 (1975).

\bibitem{Knuth_92}  D.~E. Knuth, Convolution polynomials,
   Mathematica J. {\bf 2}, 67--78 (1992),
   math.CA/9207221 at arXiv.org.

\bibitem{Labelle_80}  G. Labelle, Sur l'inversion et l'it\'eration continue
   des s\'eries formelles,
   Europ. J. Combin. {\bf 1}, 113--138 (1980).

\bibitem{Labelle_81}  G. Labelle, Une nouvelle d\'emonstration combinatoire
   des formules d'inversion de Lagrange,
   Adv. Math. {\bf 42}, 217--247 (1981).

\bibitem{Mullin_70}  R. Mullin and G.-C. Rota,
    On the foundations of combinatorial theory. III.~Theory of binomial
    enumeration,
    in {\em Graph Theory and its Applications}\/,
    edited by Bernard Harris
    (Academic Press, New York, 1970), pp.~167--213.

\bibitem{Olive_79}  G. Olive, Binomial functions and combinatorial mathematics,
   J. Math. Anal. Appl. {\bf 70}, 460--473 (1979).

\bibitem{Pemantle_13}  R. Pemantle and M.~C. Wilson,
   {\em Analytic Combinatorics in Several Variables}\/
   (Cambridge University Press, Cambridge, 2013).

\bibitem{Petreolle-Sokal_nontri}
   M. P\'etr\'eolle, A.D. Sokal and B.-X. Zhu,
   Non-triangular linear transforms preserving Hankel-total positivity,
   in preparation.


\bibitem{Roman_78}  S.~M. Roman and G.-C. Rota, The umbral calculus,
   Adv. Math. {\bf 27}, 95--188 (1978).

\bibitem{Roman_84}  S.~M. Roman, {\em The Umbral Calculus}\/
   (Academic Press, New York, 1984).

\bibitem{Rota_73}  G.-C. Rota, D. Kahaner and A. Odlyzko,
   On the foundations of combinatorial theory. VIII.~Finite operator calculus,
   J. Math. Anal. Appl. {\bf 42}, 684--760 (1973).
   [Reprinted in G.-C. Rota, {\em Finite Operator Calculus}\/
    (Academic Press, New York--London, 1975), Chapter~2.]

\bibitem{Rota_94}  G.-C. Rota and B.~D. Taylor, The classical umbral calculus,
   SIAM J. Math. Anal. {\bf 25}, 694--711 (1994).

\bibitem{Scott-Sokal_expidentities}  A.~D. Scott and A.~D. Sokal,
   Some variants of the exponential formula, with application to the
   multivariate Tutte polynomial (alias Potts model),
   S\'eminaire Lotharingien de Combinatoire {\bf 61A}, article 61Ae (2009).

\bibitem{Shapiro_91}  L.~W. Shapiro, S. Getu, W.~J. Woan and L.~C. Woodson,
   The Riordan group,
   Discrete Appl. Math. {\bf 34}, 229--239 (1991).

\bibitem{Sokal_totalpos}  A.~D. Sokal, Coefficientwise total positivity
   (via continued fractions) for some Hankel matrices of combinatorial
   polynomials, in preparation.

\bibitem{Sprugnoli_94}  R. Sprugnoli, Riordan arrays and combinatorial sums,
   Discrete Math. {\bf 132}, 267--290 (1994).

\bibitem{Stanley_99}  R.~P. Stanley, {\em Enumerative Combinatorics}\/,
      vol.~2 (Cambridge University Press, Cambridge--New York, 1999).

\bibitem{Vandermonde_1772}  A.-T. Vandermonde, M\'emoire sur des irrationnelles
   de diff\'erens ordres avec une application au cercle,
   M\'emoires de Math\'ematique et de Physique,
   Tir\'es des Registres de l'Acad\'emie Royale des Sciences (1772),
   489--498.
   
\bibitem{Wang_08}  W. Wang and T. Wang, Generalized Riordan arrays,
   Discrete Math. {\bf 308}, 6466--6500 (2008).

\bibitem{Wilf_94}  H.S. Wilf, {\em generatingfunctionology}\/,
    2nd ed. (Academic Press, San Diego--London, 1994).

\bibitem{Zeng_96}  J. Zeng, Multinomial convolution polynomials,
   Discrete Math. {\bf 160}, 219--228 (1996).

\bibitem{Zhu_2006}  Zh\={u} Sh\`{\i}ji\'e [= Chu Shih-chieh],
{\em Jade Mirror of the Four Unknowns}\/, vols.~I and~II,
translated into Modern Chinese by Guo Shuchun,
translated into English by Ch'en Tsai Hsin,
revised and supplemented by Guo Jinhai
(Liaoning Education Press, Shenyang, 2006).
Originally published in Chinese at Yangzhou, 1303.

\end{thebibliography}
\end{document}